\documentclass[12pt]{amsart}
\usepackage{amsfonts}
\usepackage{amssymb}

\newcommand{\R}{\mathbb R}

\newcommand{\X}{\mathfrak X}

\newtheorem{thm}{Theorem}[section]

\newtheorem{lem}[thm]{Lemma}
\newtheorem{prop}[thm]{Proposition}
\theoremstyle{definition}

\theoremstyle{remark}
\newtheorem{rem}[thm]{Remark}

\begin{document}

\title[HERMITIAN MANIFOLDS OF CONSTANT ANTIHOLOMORPHIC SECTIONAL CURVATURES ]
{HERMITIAN MANIFOLDS OF POINTWISE CONSTANT ANTIHOLOMORPHIC SECTIONAL CURVATURES }%
\author{Georgi Ganchev and Ognian Kassabov}%
\address{Bulgarian Academy of Sciences, Institute of Mathematics and Informatics,
Acad. G. Bonchev Str. bl. 8, 1113 Sofia, Bulgaria}%
\email{ganchev@math.bas.bg}%
\address{Higher transport School, Sofia, Bulgaria}
\email{okassabov@vtu.bg}

\begin{abstract}
In dimension greater than four, we prove that if a Hermitian non-Kaehler manifold is
of pointwise constant antiholomorphic sectional curvatures, then it is of constant
sectional curvatures.
\end{abstract}
\maketitle

 \section{Introduction}
Let $(M,g,J) \, (\dim M=2n \geq 4)$ be an almost Hermitian manifold. Any  two-plane
(section) $E$ in the  tangential space $T_pM, \; p\in M$ determines  an angle
$\theta=\angle(E,JE), \; \theta \in [0, \frac{\pi}{2}]$. Two types of planes with
respect to the angle $\theta$ are remarkable: \emph{holomorphic} sections -
characterized by the condition $\theta = 0$ or $E=JE$; \emph{antiholomorphic} sections
- characterized by the condition $\theta=\frac{\pi}{2}$ or $E\perp JE$. The latter
are also known as \emph{totally real} in view of the condition $E\perp JE$.

If $\Phi$ is the fundamental K\"ahler form of the manifold, then any antiholomrphic
section $E$ is characterized by the condition $\Phi_{|E}=0$. Because of this
characterization, these tangent planes are also known as \emph{Lagrangean}.

An almost Hermitian manifold is said to be of pointwise constant antiholomorphic sectional
curvature $\nu$ if the Riemannian sectional curvature $K(E;p)$ does not depend on the
antiholomorphic section $E$ in $T_pM, \; p\in M$, i.e. $K(E;p)=\nu(p)$ is only a function
of the point $p \in M$.

A tensor characterization for an almost Hermitian manifold of pointwise constant
antiholomorphic sectional curvature in $\dim M \geq 4$ has been found in \cite{G}.

In \cite{K} it has been proved that the antiholomorphic sectional curvature $\nu(p)$
is a constant on the manifold under the condition $\dim M>4$.

A complete classification of compact Hermitian surfaces $(\dim M =4)$ with pointwise
constant antiholomorphic sectional curvature has been given in \cite{AGI}.
Four-dimensional almost Hermitian manifolds of pointwise constant antiholomorphic
sectional curvature have been studied in \cite{S}.

In this paper we consider the class of Hermitian manifolds and prove our main
\vskip 2mm
\textbf{Theorem A}. \emph{If a Hermitian non-K\"ahler manifold with a real dimension
greater than four is of pointwise constant antiholomorphic sectional curvature, then
the manifold is of constant sectional curvature.}

\section{Preliminaries}
Let $(M,g,J) \; (\dim M =2n \geq 4$ be an almost Hermitian manifold with metric $g$
and almost complex structure $J$. The tangent space to $M$ at an arbitrary point
$p \in M$ is denoted by $T_pM$ and the algebra of all differentiable vector fields
on $M$ is denoted by ${\X}M$. The K\"ahler form $\Phi$ of the structure $(g,J)$ is
defined by the equality
$$\Phi(X,Y)=g(JX,Y), \quad X,Y \in T_pM, p \in M.$$
The Levi-Civita connection of the metric $g$ is denoted by $\nabla$ and the Riemannian
curvature tensor $R$ of type (1,3) is given by $R(X,Y)Z=\nabla_X\nabla_YZ-\nabla_Y\nabla_XZ
-\nabla_{[X,Y]}Z, \; X,Y,Z \in {\X}M.$ The corresponding curvature tensor of type (0,4)
is given by $R(X,Y,Z,U)=g(R(X,Y)Z,U)$, for all vector fields $X,Y,Z,U$.

Let $\{e_1,...,e_{2n}\}$ be an orthonormal basis at a point $p\in M$. The Ricci tensor
$\rho$ and the scalar curvature $\tau$ of the metric $g$ are determined as follows
$$\rho(X,Y)=\sum_{i=1}^{2n}R(e_i,X,Y,e_i), \quad \tau =\sum_{i=1}^{2n}\rho(e_i,e_i);
\quad X,Y \in T_pM.$$
The almost complex structure $(g,J)$ gives rise to the $*$-Ricci tensor $\rho^*$ and to
the $*$-scalar curvature $\tau^*$ defined by the formulas
$$\rho^*(X,Y)=\sum_{i=1}^{2n}R(e_i,X,JY,Je_i), \quad
\tau^*=\sum_{i=1}^{2n}\rho^*(e_i,e_i); \quad X,Y \in T_pM.$$
While the Ricci tensor is symmetric, the $*$-Ricci tensor has the property
$$\rho^*(JX,JY)=\rho(Y,X), \quad X,Y \in T_pM. \leqno(2.1)$$

The following tensor of type (0,3)
$$F(X,Y,Z)= g((\nabla_XJ)Y,Z), \quad X,Y,Z \in {\X}M $$
is closely related to the structure $(g,J)$.
This tensor satisfies the following properties
$$F(X,Y,Z)=-F(X,Z,Y), \quad F(X,JY,JZ)=-F(X,Y,Z). \leqno(2.2)$$
The well known classes of almost Hermitian manifolds have been obtained
in terms of the properties of the tensor $F$ in \cite{GH}.

In this section we consider Hermitian manifolds, which are characterized by the
following property of the tensor $F$ \cite{GH}:
$$(\nabla_{JX}J)Y=J(\nabla_XJ)Y \quad \iff \quad F(JX,Y,Z)=-F(X,JY,Z). \leqno(2.3)$$

Let $T_p^{\mathbb{C}}M$ be the complexification of the tangent space $T_pM$ at
any point $p \in M$. By ${\X}^{\mathbb{C}}M$ we denote the algebra of complex
differentiable vector fields on $M$. The complex structure $J$ generates the
standard splittings
$$T_p^{\mathbb{C}}M=T_p^{1,0}M \oplus T_p^{0,1}M, \quad
{\X}^{\mathbb{C}}M={\X}^{1,0}M\oplus{\X}^{0,1}M.$$

If $\{e_1,...,e_n; Je_1,...,Je_n\}$ is an orthonormal frame at a point $p\in M$, then
the vectors $\displaystyle{Z_{\alpha}=\frac{e_{\alpha}-iJe_{\alpha}}{2}}$ and
$Z_{\bar\alpha}=\bar Z_{\alpha}=\displaystyle{\frac{e_{\alpha}+iJe_{\alpha}}{2}};
\;\alpha =1,...,n$ form a basis for $T_p^{1,0}M$ and $T_p^{0,1}M$, respectively.
Further, we call these bases $\{Z_{\alpha};Z_{\bar\alpha}\}\; \alpha=1,...,n$
\emph{special complex bases}.

For an arbitrary tensor T we denote $T_{\alpha ...}=T(Z_{\alpha} ...)$ and
$T_{\bar \alpha ...}=T(Z_{\bar\alpha}...).$

In what follows, the summation convention is assumed and Greek indices
$\alpha, \beta, \gamma,...$ run from $1$ to $n$.

It follows that the components of the metric tensor with respect to a special complex
basis satisfy the conditions
$$g_{\alpha\beta}=0 \, ({\rm for \, all} \, \alpha, \beta); \quad  g_{\alpha\bar\beta}=0,\;
({\rm for} \, \alpha\neq \beta); \quad g_{\alpha\bar\alpha}=
\frac{1}{2} \ ({\rm for \, all} \, \alpha).$$

We have the following
\begin{lem}\label{L:2.1}
Let $(M,g,J)$ be a Hermitian manifold. If $(\nabla_{\bar Z}J)Z=0$ for an arbitrary
$Z\in T_p^{1,0}M$, then $\nabla J=0$ at the point $p$.
\end{lem}

\emph{Proof:} Indeed, the condition $(\nabla_{\bar Z}J)Z=0$ implies that
$(\nabla_XJ)X=0$ for all $X\in T_pM$. Therefore $M$ satisfies the condition
characterizing a nearly K\"ahler manifold at $p$. Since $M$ is Hermitian, then $M$
is K\"ahlerian, i.e. $\nabla J=0$ at $p$ \cite{GH} .
\hfill{\qed}
\vskip 2mm
Now, let $(M,g,J)$ be a Hermitian manifold with pointwise constant antiholomorphic
sectional curvature. This means that for any orthonormal antiholomorphic frame
$\{X,Y\}$, \; $(g(X,X)=g(Y,Y)=1, \, g(X,Y)=g(X,JY)=0)$ at an arbitrary point $p \in M$
the sectional curvature $R(X,Y,Y,X)$ does not depend on the antiholomorphic section
${\rm span}\{X,Y\}$, i.e. $R(X,Y,Y,X)$ is only a function of the point $p$.
We denote this function by $\nu(p)$.

Let $Q(X,Y)$ be a tensor on $M$ having the symmetry (2.1), i.e.
$$Q(JX,JY)=Q(Y,X).\leqno(2.4)$$

The following tensor construction $\Psi(Q)$ is relevant to the considerations in this paper:
$$\begin{array}{ll}
\Psi(Q)(X,Y,Z,U)=&\;g(Y,JZ)\,Q(X,JU)-g(X,JZ)\,Q(Y,JU)-2\,g(X,JY)\,Q(Z,JU)\\
[3mm]
&+g(X,JU)\,Q(Y,JZ)-g(Y,JU)\,Q(X,JZ)-2\,g(Z,JU)\,Q(X,JY).\end{array}$$

We also recall the basic invariant tensors $\pi_1$ and $\pi_2$ only formed by the
fundamental tensors $g$ and $\Phi$:
$$\pi_1(X,Y,Z,U)=g(Y,Z)g(X,U)-g(X,Z)g(Y,U),$$
$$\pi_2(X,Y,Z,U)=g(Y,JZ)\,g(X,JU)-g(X,JZ)\,g(Y,JU)-2\,g(X,JY)\,g(Z,JU).$$

The first author has proved the following tensor characterization for an
almost Hermitian manifold of pointwise constant antiholomorphic sectional curvatures.
\vskip 2mm
\textbf{Theorem.\cite{G}} \emph{An almost Hermitian manifold with $\dim M = 2n\geq 4$ is of pointwise
constant antiholomorphic sectional curvature $\nu(p)$ if and only if its curvature
tensor satisfies the identity}
$$R-\frac{1}{2(n+1)}\,\Psi(\rho^*)+\frac{\tau^*}{2(n+1)(2n+1)}\,\pi_2=
\nu\,\left(\pi_1-\frac{1}{2n+1}\,\pi_2\right).\leqno(2.5)$$
\vskip 2mm
We introduce the tensor
$$Q=\frac{1}{2(n+1)}\,\rho^*-\frac{\tau^*+2(n+1)\nu}{4(n+1)(2n+1)}\,g\, ,$$
which in view of (2.1) has the property (2.4). Then the condition (2.5) can be written as follows:
$$R=\Psi(Q)+\nu \, \pi_1.\leqno(2.6)$$

The second author has proved in \cite{K} that in $\dim M\geq 6$ the function $\nu(p)$ in (2.5)
is constant. Thus, we shall speak about almost Hermitian manifolds of constant antiholomorphic
sectional curvature instead of "pointwise constant" antiholomorphic sectional curvature.

\section{Proof of Theorem A}

In this section we prove Theorem A on the base of the following statement.
\begin{prop}\label{P:3.1}
Let $(M,g,J) \; (\dim_{\mathbb C} M \geq 3)$ be a Hermitian  manifold of constant
antiholomorphic sectional curvature. Then any non-K\"ahler point of $M$ has a neighborhood
in which $(M,g,J)$ is of constant sectional curvature.
\end{prop}
\textbf{Proof:} Let $p_0 $ be a point in $M$ with $F\neq 0$ at $p_0$. We consider a
neighborhood $U$ of $p_0$, such that $F\neq 0$ at any point of $U$. We shall prove that
$(M,g,J)$ is of constant sectional curvature in $U$.

For any $p\in U$, we consider a special complex basis $\{Z_{\alpha}, Z_{\bar\alpha}\}
\; \alpha=1,...,n$ at the point $p$.

The property (2.4) of the tensor $Q$ implies that
$$Q_{\alpha \bar \beta}=Q_{\bar \beta \alpha}, \qquad
Q_{\alpha\beta}=-Q_{\beta\alpha}.\leqno(3.1)$$

Taking into account the property (2.3) of the covariant derivative of the complex structure
and the symmetry (2.4) of the tensor $Q$, we compute
$$(\nabla_XQ)(JY,JZ)=(\nabla_XQ)(Z,Y)-Q((\nabla_XJ)Y,JZ)-Q(JY,(\nabla_XJ)Z)\leqno(3.2)$$
for arbitrary $X,Y,Z \in {\X}U$.

Since the tensor $F$ has the symmetries (2.2) and (2.3), then its essential components
(those which may not be zero) with respect to a special complex basis
$\{Z_{\alpha}, \; Z_{\bar \alpha}\}$ are only $F_{\bar\alpha\beta\gamma}$ and their conjugates.
These components satisfy the condition $F_{\bar\alpha\beta\gamma}= -F_{\bar\alpha\gamma\beta}$.
These properties of the tensor $F$ can be expressed in terms of the covariant derivative
$(\nabla_XJ)Y$ as follows
$$\nabla_{\alpha}J_{\beta}^{\gamma}= \nabla_{\bar\alpha}J_{\beta}^{\gamma}=
\nabla_{\alpha}J_{\beta}^{\bar\gamma}=0.\leqno(3.3)$$

The equalities (3.2) and (3.3) imply that
$$\nabla_{\alpha}Q_{\bar\gamma\beta}=\nabla_{\alpha}Q_{\beta\bar\gamma}+
i\nabla_{\alpha}J_{\bar\gamma}^{\sigma}Q_{\beta\sigma};\leqno(3.4)$$
$$\nabla_{\alpha}Q_{\beta\gamma}=-\nabla_{\alpha}Q_{\gamma\beta}, \quad ({\rm especially}
\quad \nabla_{\alpha}Q_{\beta\beta}=0); \leqno(3.5)$$
$$\nabla_{\bar\alpha}Q_{\beta\beta}=
i\nabla_{\bar\alpha}J_{\beta}^{\bar\sigma}Q_{\bar \sigma\beta}.\leqno(3.6)$$

First we prove the following statement
\begin{lem}\label{L:3.1}
Let $Z,W \in T_p^{1,0}M$ and $g(Z,\bar W)=0$. If $F(\bar Z,Z, W)\neq 0$, then $Q(Z,W)=0$.
\end{lem}
\textbf{Proof:} Since $g(Z,\bar W)=0$, then we can find a special complex basis
$\{Z_{\alpha}, Z_{\bar\alpha}\} \; \alpha=1,...,n$ such that the vectors $Z$ and $W$
are collinear with $Z_{\alpha}$ and $Z_{\beta}$, respectively, for some $\alpha \neq \beta$.

Applying the Bianchi identity for the curvature tensor $R$ in the form
$$\nabla_{\alpha}R_{\beta\gamma\beta\bar\gamma}+\nabla_{\beta}R_{\gamma\alpha\beta\bar\gamma}
+\nabla_{\gamma}R_{\alpha\beta\beta\bar\gamma}=0,$$
we find
$$\nabla_{\beta}Q_{\alpha\beta}=0.\leqno(3.7)$$

Further we apply the Bianchi identity in the form
$$\nabla_{\bar\alpha}R_{\alpha\beta\alpha\beta}+\nabla_{\alpha}R_{\beta\bar\alpha\alpha\beta}
+\nabla_{\beta}R_{\bar\alpha\alpha\alpha\beta}=0$$
and taking into account (3.7), we obtain
$$F_{\bar\alpha\alpha\beta}Q_{\alpha\beta}=0. \leqno(3.8)$$

Under the conditions of the lemma we have $F_{\bar\alpha\alpha\beta}\neq 0$. Then it
follows from (3.8) that $Q_{\alpha\beta}=0$.
\hfill{\qed}
\vskip 2mm
Next we prove
\begin{lem}\label{L:3.2}
The tensor $Q$ is symmetric at any point $p \in U$.
\end{lem}
\textbf{Proof:} Since the tensor $F\neq 0$ at the point $p$, then because of the Lemma
\ref{L:2.1} there exist indices $\alpha\neq\beta$ so that $F_{\bar\alpha\alpha\beta}\neq 0$.
Applying Lemma \ref{L:3.1}, it follows that $Q_{\alpha\beta}=0$.

Let $\gamma\neq \alpha, \beta$. Since $F_{\bar\alpha\alpha\beta}\neq 0$, then the
complex function $w(t)=F(Z_{\bar\alpha},Z_{\alpha}, Z_{\beta}+tZ_{\gamma})\neq 0$ for all
sufficiently small $t\in {\R}$. It follows from Lemma \ref{L:3.1} that $Q(Z_{\alpha},
Z_{\beta}+tZ_{\gamma})=0$. Hence, $Q(Z_{\alpha}, Z_{\gamma})=0$, i.e. $Q_{\alpha\gamma}=0$.

Similarly, the inequality $F(Z_{\alpha}+tZ_{\gamma},Z_{\alpha}+tZ_{\gamma},Z_{\beta})\neq 0$
for all sufficiently small real $t$ and Lemma 3.2 imply that
$Q(Z_{\alpha}+tZ_{\gamma},Z_{\beta})=0$.
Hence, $Q_{\gamma\beta}=0$. So far, we obtained
$$Q_{\alpha\beta}=Q_{\alpha\gamma}=Q_{\beta\gamma}=0.$$

In $\dim M > 6$, let $\delta\neq \alpha, \beta, \gamma$. As in the above, we find
$$Q_{\alpha\delta}=Q_{\alpha\delta}=Q_{\beta\delta}=0.$$

On the other hand, the inequality $F(Z_{\bar\alpha}+tZ_{\bar\gamma}, Z_{\alpha}+tZ_{\gamma},
Z_{\beta}+tZ_{\delta}\neq 0)$, which is valid for sufficiently small real $t$, implies that
$Q_{\gamma\delta}=0$.

Thus we obtained $Q_{\lambda\mu}=0$ for all $\lambda, \mu=1,...,n$, which proves
the assertion.
\hfill{\qed}
\vskip 2mm
Finally, we shall prove that the tensor $Q$ is proportional to the metric tensor $g$ in $U$.

For that purpose it is sufficient to prove that
$$Q_{\lambda\bar\mu}=0, \leqno(3.9)$$
for all different indices $\lambda$ and $\mu$.

We consider two cases for the tensor $F\neq 0$:

1) There exist three different indices $\alpha, \beta, \gamma$, such that
$F_{\bar\gamma\alpha\beta} \ne 0$;

2) $F_{\bar\gamma\alpha\beta} = 0$ for all different indices $\alpha, \beta, \gamma$
with respect to any special complex basis.
\vskip 2mm
The case 1). Applying the second Bianchi identity in the form
$$ \nabla_{\alpha} R_{\beta\bar\gamma\beta\bar\gamma} +
   \nabla_{\beta} R_{\bar\gamma\alpha\beta\bar\gamma} +
   \nabla_{\bar\gamma} R_{\alpha\beta\beta\bar\gamma}= 0, $$
we get the equality $ F_{\bar\gamma\alpha\beta}Q_{\beta\bar\gamma} =0 ,$
which implies that $Q_{\beta\bar\gamma} =0$.

Now, arguments similar to those in Lemma \ref{L:3.2} show (3.9).
\vskip 2mm
The case 2). According to Lemma \ref{L:2.1} there exist two different indices
$\alpha$ and $\beta$ such that $F_{\bar\alpha\alpha\beta}\neq 0$.
Applying the second Bianchi identity in the form
$$\nabla_{\alpha} R_{\bar\gamma\beta\beta\bar\alpha} +
  \nabla_{\bar\gamma} R_{\beta\alpha\beta\bar\alpha} +
  \nabla_{\beta} R_{\alpha\bar\gamma\beta\bar\alpha}= 0, $$
and taking into account the equalities $Q_{\alpha\beta}=0$,
$ F_{\bar\gamma\alpha\beta}=0 $, we find
$$  - \nabla_{\bar\gamma} Q_{\beta\beta}
    + iQ(Z_\beta,(\nabla_{\bar\gamma} J)Z_\beta)
    + \nabla_{\beta} Q_{\bar\gamma\beta}=0. $$

The last equality in view of (3.6) implies
$$  \nabla_{\beta} Q_{\bar\gamma\beta} =0. \leqno (3.10)$$

Applying the second Bianchi identity in the form
$$ \nabla_{\bar\alpha} R_{\alpha\beta\beta\bar\gamma} +
   \nabla_{\alpha} R_{\beta\bar\alpha\beta\bar\gamma} +
   \nabla_{\beta} R_{\bar\alpha\alpha\beta\bar\gamma}= 0, $$
we find
$$  3iF_{\bar\alpha\alpha\beta}Q_{\beta\bar\gamma}
   +2\nabla_{\beta} Q_{\beta\bar\gamma}
   +2iQ(Z_\beta,(\nabla_\beta J)Z_{\bar\gamma})=0, $$
which together with (3.4) gives
$$  3iF_{\bar\alpha\alpha\beta}Q_{\beta\bar\gamma}
   +2\nabla_{\beta} Q_{\bar\gamma\beta} =0\ . $$

The last equality and (3.10) imply that
$$  F_{\bar\alpha\alpha\beta}Q_{\beta\bar\gamma}=0 \ . $$

Hence, $Q_{\beta\bar\gamma}=0 $.
Applying again the scheme of the proof of Lemma \ref{L:3.2}, we obtain the
conditions (3.9).

Thus, in both cases 1) and 2), we obtained the conditions (3.9),
which are equivalent to the identity
$$Q(X,Y)=0, \; {\rm whenever} \; X,Y\in T_pM, \; g(X,Y)=0.\leqno(3.11)$$
Applying standard arguments for the symmetric tensor $Q(X,Y)$, we obtain
that the tensor $Q$ is proportional to the metric tensor $g$, i.e.
$$Q=\frac{{\rm tr}\,Q}{2n}\,g, \quad {\rm tr}\,Q=\frac{\tau^*-2\,n\nu}{2(2n+1)}\,.$$
Hence
$$R=\nu \,\pi_1+\frac{{\rm tr} \, Q}{n} \,\pi_2.$$

Further we use the following statement
\vskip 2mm
\textbf{Theorem. \cite{TV}} \emph{Let $M$ be a connected almost Hermitian manifold
with real dimension $2n\geq 6$ and Riemannian curvature tensor of the following form:
$$R=f\,\pi_1+h\,\pi_2,$$
where $f$ and $h$ are ${\mathcal C}^{\infty}$ functions on $M$ such that $h$ is not
identical zero. Then $M$ is a complex space form (i.e. a K\"ahler manifold with
constant holomorphic sectional curvature).}
\vskip 2mm

Applying the above mentioned theorem, we obtain that the function ${\rm tr} \, Q=0$,
i.e. $\tau^*-2n\nu=0$. Hence $M$ is of constant sectional curvature $\nu$ in $U$.
\hfill{\qed}
\vskip 2mm
\begin{rem}
If the curvature tensor of an almost Hermitian manifold has the form $R=\nu \, \pi_1$,
then
$$\nu=\frac{\tau}{2n(2n-1)}=\frac{\tau^*}{2n}$$
and $\tau=(2n-1)\tau^*$.
\end{rem}
To complete the proof of Theorem A, denote by $H$ the set of points in $M$, in which
$F\neq 0$. Then $H$ is a non-empty open set of $M$ and according to Proposition 3.1
$M$ is of constant sectional curvature $\nu$, i.e. $R=\nu\,\pi_1$ in $H$.

Let $K=\{ p \in M \ : \ R \ne \nu\pi_1 \}$. Then $K$ is also open and we have
$ H\cap K = \phi$, $ H\cup K = M$. Since $M$ is connected and $H$ is nonempty, then
$K$ is empty, which completes the proof of Theorem A.

\hfill{\qed}

\end{document}